\documentclass[10pt]{article}
\newtheorem{df}{Definition}
\newtheorem{lem}[df]{Lemma}

\newtheorem{ex}[df]{Example}

\title{Fixes of permutations acting on monotone Boolean functions.}
\author{
\Large Andrzej Szepietowski\\
\small Institute of Informatics, \\[-0.8ex] 
\small Faculty of Mathematics, Physics and Informatics, \\[-0.8ex]
\small University of Gda\'{n}sk, \\[-0.8ex]
\small 80-308 Gda\'{n}sk, Poland\\[-0.8ex]
\small {\tt matszp@inf.ug.edu.pl}\\
}

\date{}

\begin{document}

\maketitle 

\begin{abstract}
We present a few algorithms and methods to count fixes of permutations acting on monotone Boolean functions. Some of these  methods was used by Pawelski~\cite{P} to compute the number of inequivalent monotone Boolean functions with 8 variables.
\end{abstract}


\section{Introduction}
Let $B$ denote the set $\{0,1\}$ and $B^n$ the
set of $n$-element sequences of $B$.
A Boolean function with $n$ variables is any function
from $B^n$ into $B$. There are $2^n$ elements in $B^n$
and $2^{2^n}$ Boolean functions with $n$ variables.
There is the order relation in $B$ (namely: $0\le 0$, $\;0\le1$, $\;1\le1$)
and the partial order in $B^n$:
for any two elements: 
$x=(x_1,\dots,x_n)$, $\;y=(y_1,\dots,y_n)$ in $B^n$,
$$x\le y\quad\hbox{if and only if}\quad x_i\le y_i
\quad \hbox{for all $1\le i\le n$}.$$
The function $h:B^n\to B$ is monotone if
$$x\le y \Rightarrow h(x)\le h(y).$$
Let $D_n$ denote the set of monotone functions with $n$ variables
and let $d_n$ denote $|D_n|$.
Known values of $d_n$, for $n=0,\dots,8$ are presented in the table at the end of this paper.
The values $d_n$ for $n\le 4$ was published by Dedekind~\cite{D}, Church~\cite{Ch1,Ch2} gave the values $d_5$ and
$d_7$, Ward~\cite{Wa} the value $d_6$, and the last known value $d_8$ was  published by Wiedemann~\cite{W}.
Dedekind numbers was also considered in~\cite{B,C,CW,FMSS, SY}.

We have  the partial order in $D_n$ defined by:
$$g\le h\quad\hbox{if and only if}\quad g(x)\le h(x)\quad
\hbox{for all } x\in B^n.$$
We shall represent the  elements of $D_n$ as strings of bits of
length $2^n$. Two elements of $D_0$ will be represented as 0 and 1;
any element $g\in D_1$ can  be represented as 
the concatenation $g(0)*g(1)$, where $g(0), g(1)\in D_0$ and $g(0)\le g(1)$. Hence, $D_1=\{00, 01, 11\}$. 
Each element of $g\in D_2$ is the concatenation (string) of four bits: $g(00)*g(10)*g(01)*g(11)$ which can be  represented as
a concatenation $g_0*g_1$, where $g_0, g_1\in D_1$ and $g_0\le g_1$.
Hence, $D_2=\{0000, 0001, 0011, 0101, 0111, 1111\}$.
Similarly any element of $g\in D_n$ can be represented as
a concatenation $g_0*g_1$, where $g_0, g_1\in D_{n-1}$ and $g_0\le g_1$.

Let $S_n$ denote the set of permutations on $\{1, ... n\}$.
Every permutation $\pi\in S_n$ defines the permutation on $B^n$ by
$\pi(x)=x\circ\pi$
(we treat each element $x\in B^n$ as a function
$x:\{1, ... ,n\}\to \{0,1\}$).
Note that $x\le y$ if and only if $\pi(x)\le\pi(y)$.
This permutation on $B^n$ generates the permutation
on $D_n$. Namely, by $\pi(g)=g\circ \pi$.
Note that $\pi(g)$ is monotone if $g$ is monotone.
Two functions $f, g\in D_n$ are equivalent if there is a 
permutation $\pi\in S_n$ such that $f=\pi(g)$.
By $R_n$ we denote the set of equivalent classes and by $r_n$
we denote the number of the equivalence classes.
Known values of $r_n$ (for $n\le 8$) are given in the table at the end of this paper.
The number of the equivalence classes can be computed by Burnside's Lemma, see~\cite{M}. Namely,
$$
r_n={1\over n!}\sum_{\pi\in S_n}|Fix(\pi,D_n)|,
$$ 
where $Fix(\pi,D_n)$ is the set of fixes of the permutation $\pi$ acting on $D_n$.
A function $f\in D_n$ is a fix of $\pi$ if $\pi(f)=f$. It is obvious that for two permutations
$\pi$ and  $\rho$ of the same cycle index, $|Fix(\pi,D_n)|=| Fix(\rho,D_n)|$. Hence,
$$
r_n={1\over|n!|}\sum_{i=1}^r \mu_i\cdot|Fix(\pi_i,D_n)|
$$ 
where $\pi_i$ is a permutation of the cycle index $i$ and  $\mu_i$ is the number of premutations of index $i$.
 
 In 1985 and 1986 Chuchang and Shoben~\cite{CS1,CS2} used Burnside Lemma in order to compute $r_n$ for all $n\le 7$.
 Recently, Pawelski~\cite{P} presented $r_8$. 
  
In this paper we present a few algorithms and methods to count fixes
of permutations acting on $D_n$. Some of these methods were used by Pawelski in~\cite{P}.

\section{Posets}
A poset (partially ordered set) $(S,\le)$ consists of a set $S$ together with a binary relation (partial order) $\le$ which is reflexive, transitive and antisymmetric.
For two posets $(S,\le)$ and $(T,\le)$ by $S\times T$ we denote the cartesian product with the order defined by: $(a,b)\le(c,d) $ iff $a\le c $ and $b\le d$.
For two disjoint posets $(S,\le)$ and $(T,\le)$ by $S+ T$ we denote the disjoint union (sum)  with the order defined by:
$$x\le y\quad\hbox{iff}\quad (x,y\in S \quad\hbox{and}\quad x\le y)\quad
\hbox{or}\quad(x,y\in T\quad\hbox{ and }\quad x\le y).$$
Given two posets $(S,\le)$ and $(T,\le)$ a function $f:S \to T$ is called monotone, if for all $x,y\in S$,  $x \le y$ implies $f(x)\le f(y)$. 
By $T^S$ we denote  the poset of all monotone functions from $S$ to $T$ with the partial order  defined by:
$$f\le g\quad\hbox{ if and only if}\quad f(x)\le g(x)\hbox{ for all }x\in S.$$

\begin{lem}\label{L1} For three posets $R,S,T$:

(1) If $S$ and $T$ are disjoint, then the poset $R^{S+T}$ is isomorphic to $R^S\times R^T$.

(2) The poset $R^{S\times T}$ is isomorphic to $(R^S)^T$ and to $(R^T)^S$.
\end{lem}

In this paper we shall deal with the following posets:
\begin{itemize}
\item $A_n=\{1,\dots,n\}$ antichain,  where no two distinct elements are related.
\item $B=\{0,1\}$ two bits ordered by  $0\le 0$, $0\le1$,  $1<1$.
\item $B^n=B^{A_n}$ which is isomorphic to:
\begin{itemize}
\item the poset of all subsets of $\{1,\dots,n\}$ ordered by the inclusion,
\item the poset of all $n$-strings of bits, where $(x_1,\dots,x_n)\le(y_1,\dots,y_n)$ iff $x_i\le y_i$ for all $i$. 
\end{itemize}
\item $D_n=B^{B^n}$ the poset of all monotone Boolean functions of $n$ variables.
\item Path (or chain) $P_n=\{p_1<\dots<p_n\}$. Note that $B^{P_n}=P_{n+1}$.
\end{itemize}

\section{Arrays}
Let $M(S)$ denote the array of the poset $S$. For $i,j\in S$, $M(S)[i,j]=1$ if $i\le j$, and  $M(S)[i,j]=0$ otherwise.
For example, for the poset
$D_1=\{00<01<11\}$ its array
$$
M(D_1)=\left(
\begin{array}{ccc}
1 & 1& 1\\
0 & 1& 1\\
0 & 0&1
\end{array}
\right)
$$
The poset $D_1$ is equal (isomorphic) to the poset of the path $P_3=\{a<b<c\}$. 

The elements of $M(S)$ describe monotone functions from the poset $B=\{0,1\}$ to $S$.
If $M(S)[i,j]=1$ then there exists the monotone function with $f(0)=i$ and $f(1)=j$.
Thus, if we add the elements of $M(S)$ we obtain $|S^B|$ - the number of monotone functions
from $B$ to $S$. For example
$$Sum(M(D_1))=6=|D_1^B|=|(B^B)^B|= |B^{B\times B}|=|B^{B^2}| =|D_2|=d_2,$$
where $Sum(M(D_1))$ denotes the sum of all elements of the array $M(D_1)$.

Consider the product $M(S)^2=M(S)\times M(S)$. Then
$M(S)^2[i,j]=|\{k\in S : i\le k\le j\}|$ which is the number of elements in the interval 
$[i,j]\subset S$.
Moreover the elements of $M(S)^2$ are connected to monotone  functions from the path $P_3=\{a<b<c\}$ to $S$. Indeed, $M(S)^2[i,j]$ is equal to the number of monotone functions with $f(a)=i$ and $f(c)=j$. Hence, $Sum(M(S)^2)=|S^{P_3}|$.
For example,
$$
M(D_1)^2=\left(
\begin{array}{ccc}
1 & 2& 3\\
0 & 1& 2\\
0 & 0&1
\end{array}
\right)
$$

$M(D_1)^2[1,3]=3$ is equal to the number of elements in the interval 
$[00,11]=\{00,01,11\}$. Furthermore, 
$Sum(M(D_1)^2)=10$ is equal to $|D_1^{P_3}|$ - the number of monotone functions from $P_3$ to  $D_1$, and to  $|(B^{B})^{P_3}|$, and to $|B^{B\times P_3}|$.

Consider monotone functions from the cube $B^2=(00,01,10,11)$ to $S$, then  $(M(S)^2[i,j])^2=M(S)^2[i,j]\cdot M(S)^2[i,j]$ is equal to the number of monotone functions
with $f(00)=i$ and $f(11)=j$. Hence, $SumSq(M(S)^2)=|S^{B^2}|,$
where $SumSq(M(S)^2)$ denotes the sum of squares of all elements of the array $M(S)^2$.
For example $SumSq(M(D_1)^2)=20$ is equal to $|D_1^{B^2}|=|(B^B)^{B^2}|=|B^{B\times B^2}|= |B^{B^3}|=d_3$.

Consider the product $M(S)^3=M(S)\times M(S)\times M(S)$. The elements of $M(S)^3$ are connected  to
monotone functions from the path $P_4=(a<b<c<d)$ to $S$. Indeed,   $M(S)^3[i,j]$ is equal to the number of monotone functions
with $f(a)=i$ and $f(d)=j$.
Hence, $Sum(M(S)^3)=|S^{P_4}|.$
For example, 
$$
M(D_1)^3=\left(
\begin{array}{ccc}
1 & 3& 6\\
0 & 1& 3\\
0 & 0&1
\end{array}
\right)
$$
$Sum(M(D_1)^3)=15$ is equal to $|D_1^{P_4}|$ - the number of monotone functions from $P_4$ to  $D_1$, and to  $|(B^{B})^{P_4}|$, and to $|B^{B\times P_4}|$.

\section{Symmetries in $B^n$ and in $D_n$}
Let $S_n$ denote the set of permutations on $\{1,\dots,n\}$.
Every permutation $\pi\in S_n$ defines the permutation on $B^n$ by
$\pi(x)=x\circ\pi$.
Here we treat elements $x\in B^n$ as functions $x:\{1,\dots ,n\}\to \{0,1\}.$
Note that $x\le y$ if and only if $\pi(x)\le\pi(y)$.

\begin{ex} Consider two permutations $\pi_1=(1,2)$ and $\pi_2=(1,2,3)$
acting on $B^3$

\begin{tabular}[h]{|c||c|c|c|c|c|c|c|c|}

\hline
$x$       &000 &100 &010 &110 &001 &101 &011 &111\\
\hline
$\pi_1(x)$&000 &010 &100 &110 &001 &011 &101 &111\\ 
\hline
$\pi_2(x)$&000 &001 &100 &101 &010 &011 &110 &111\\ 
\hline
\end{tabular}
\end{ex}
Each permutation $\pi$ acting on $B^n$ generates the permutation
on $D_n=B^{B^n}$. Namely: $\pi(g)=g\circ \pi$.
Note that $\pi(g)$ is monotone if $g$ is monotone.
Two functions $f, g\in D_n$ are equivalent if there is a 
permutation $\pi\in S_n$ such that $f=\pi(g)$.
By $R_n$ we denote the set of equivalent classes.
The number of the equivalence classes denoted by $r_n$ can be computed by Burnside's Lemma, see~\cite{M}.
Known values of $r_n$ (for $n\le 8$) are given in the table at the end of this paper.
By the lemma the number of equivalence classes in $D_n$ is equal to:
$$
r_n={1\over n!}\sum_{\pi\in S_n}|Fix(\pi,D_n)|
$$ 
where $Fix(\pi,D_n)$ is the set of fixes of the permutation $\pi$ acting on $D_n$.
A function $f\in D_n$ is a fix of $\pi$ if $\pi(f)=f$. 
It is obvious that for two permutations
$\pi$ and  $\rho$ of the same cycle index, $|Fix(\pi,D_n)|=| Fix(\rho,D_n)|$. Hence,
$$
r_n={1\over|n!|}\sum_{i=1}^r \mu_i\cdot|Fix(\pi_i,D_n)|
$$ 
where $\pi_i$ is a permutation of cycle index $i$ and $\mu_i$ is the number of premutations of index $i$.

In this paper we present a few algorithms and methods to count fixes
of permutations acting on $D_n$. 

Consider a permutation $\pi\in S_n$ and suppose that $\pi$ when acting 
on $B^n$ is a product of cycles 
$\pi=C_1\circ\cdots\circ C_r$, then a monotone function $f:B^n\to B$ is a fix of $\pi$ if and only if $f$ is constant on every cycle $C_i$. 
Let $Cycl(\pi,B^n)$ denote the set of cycles $\{C_1,\dots, C_r\}$,
and let $\le$ be the partial order defined in the following way: $C_i\le C_j$ if and only if there exist
$x\in C_i$ and $y\in C_j$ such that $x\le y$ (we identify each cycle with the set of its elements).  Hence, the poset $Fix(\pi, D_n)$ is isomorphic
to the poset $B^{Cycl(\pi,B^n)}$ of monotone functions from $Cycl(\pi,B^n)$ to $B=\{0,1\}$ and we can represent fixes in $Fix(\pi,D_n)$ as sequences of bits of length $|Cycl(\pi,B^n)|$.

For the identity permutation $e$, each $x\in B^n$ forms  a cycle of length 1, hence
$Cycl(e,B^n)=B^n$ and $Fix(e,D_n)=D_n$.     
For $n=1$, we have $D_1=R_1$ and $r_1=d_1=3$.
For $n=2$, we have two permutations:
the identity $e$ with $Fix(e,D_2)=D_2$ and the inversion $(1,2)$ with
three cycles $C_1=(00)$, $C_2=(10,01)$, and $C_3=(11)$ which form the
path $P_3=\{C_1<C_2<C_3\}$. There are four monotone function from the path $P_3$ to $B$ and four fixes in $Fix((1,2),D_2)$. By Burnside's Lemma we have:
$$r_2={1\over2}(|Fix(e,D_2)|+|Fix((1,2),D_2)|)={1\over2}(6+4)=5.$$
Indeed, there are five equivalent classes in $D_2$:
$$R_2=\{\{0000\}, \{0001\}, \{0101,0011\},\{0111\},\{1111\}\}.$$


\section{Partition $A_n=A_k+A_m$}
Consider a partition of the antichain  $A_n=\{1,\dots,n\}$ into two disjoint antichains 
$A_k=\{1,\dots,k\}$ and $A_m=\{k+1,\dots n\}$, where $n=k+m$; and two permutations:
one $\pi$ acting on $A_k$ and $\rho$ acting on $A_m$. The cube $B^n$ is isomorphic
to the cartesian product $B^n=B^k\times B^m$. 
Suppose that we have two cycles, one $C_r$ of $\pi$ acting  on $B^k$ and the other $C_s$ of $\rho$ accting on $B^m$. Note that if the lengths of the two cycles are coprime, then the product
$C_r\times C_s$ is a cycle of $\pi\circ \rho$ acting on $B^n$.
Furthermore, if each cycle of $\pi$ has the length which is coprime with the length of
every cycle of $\rho$ then
$$Cycl(\pi\circ\rho, B^n)=Cycl(\pi, B^k)\times Cycl(\rho, B^m)$$
and
$$Fix(\pi\circ\rho,D_n)=Fix(\pi,D_k)^{Cycl(\rho, B^m)}=
Fix(\rho,D_m)^{Cycl(\pi,B^k)}.$$
\begin{lem}
Consider a partition of the antichain  $A_n=\{1,\dots,n\}$ into two disjoint antichains $A_k=\{1,\dots,k\}$ and $A_m=\{k+1,\dots n\}$, where $n=k+m$; and two permutations: one $\pi$ acting on $A_k$ and $\rho$ acting on $A_m$. Suppose that each cycle of $\pi$ has the length which is coprime with the length of
every cycle of $\rho$ then
$$Fix(\pi\circ\rho,D_n)=Fix(\pi,D_k)^{Cycl(\rho,B^m)}=
Fix(\rho,D_m)^{Cycl(\pi,B^k)}.$$
\end{lem}

\begin{ex} Consider partition $A_5=A_3+A_2$, with $A_3=\{1,2,3\}$ and $A_2=\{4,5\}$,
and two permutations $\pi=(1,2,3)$ and $\rho=(4,5)$. Furthermore, consider two cycles:
$C_r=(100,001,010\}$ of $\pi$ acting on $B^{\{1,2,3\}}$ and $C_s=\{10,01\}$ of $\rho$ acting on $B^{\{4,5\}}$. Their cartesian product
$$C_r\times C_s=\{10010,00101,01010,10001,00110,01001\}$$ forms one cycle of $\pi\circ\rho$ acting on $B^5=B^3\times B^2$.  
Cycles of $(1,2,3)$ acting on $B^3$ are of length 1 or 3 and they form the chain $P_4$. Cycles of $(4,5)$ acting on $B^2$ are of length 1 or 2 and they form the chain $P_3$.
Hence the set of cycles of $\pi\circ\rho$ acting on $B^5$ is the cartesian product
$P_4\times P_3$ and 
$$Fix((123)(45),D_5)=|B^{P_4\times P_3}|=|P_5^{P_3}|=|P_4^{P_4}|. $$
In order to compute $Fix((123)(45),D_5)$, let us consider the arrays:

$$
M(P_4)=
\left(
\begin{array}{cccc}
1& 1& 1& 1 \\
0& 1& 1& 1 \\
0& 0& 1& 1 \\
0& 0& 0& 1 
\end{array}
\right)
$$
$$
M(P_4)^3=
\left(
\begin{array}{cccc}
1& 3& 6& 10 \\
0& 1& 3& 6 \\
0& 0& 1& 3 \\
0& 0& 0& 1 
\end{array}
\right)
$$
$Fix((123)(45),D_5)=|P_4^{P_4}|=Sum(M(P_4)^3)=35$.
We can also use the arrays:

$$
M(P_5)=
\left(
\begin{array}{ccccc}
1& 1& 1& 1&1 \\
0& 1& 1& 1&1 \\
0& 0& 1& 1&1\\
0& 0& 0& 1&1\\ 
0& 0& 0& 0&1
\end{array}
\right)
$$
$$
M(P_5)^2=
\left(
\begin{array}{ccccc}
1& 2& 3& 4&5 \\
0& 1& 2& 3&4 \\
0& 0& 1& 2&3\\
0& 0& 0& 1&3\\ 
0& 0& 0& 0&1
\end{array}
\right)
$$
$Fix((123)(45),D_5)=|P_5^{P_3}|=Sum(M(P_5)^2)=35$.
\end{ex}

Consider a permutation $\pi$ acting on $\{1,\dots, k\}$ and on $B^k$.
Then we can say that $\pi$ also acts on $\{1,\dots, n\}$ and on $B^n$,
by identifying $\pi$ with $\pi\circ e$. Every cycle in $Cycl(e,B^m)$ has length 1,  $Cycl(e,B^m)=B^m$. Hence, 
$$Cycl(\pi\circ e, B^n)=Cycl(\pi, B^k)\times Cycl(e, B^m)=Cycl(\pi, B^k)\times  B^m$$
and
$$Fix(\pi,D_n)= Fix(\pi\circ e, D_n)=B^{Cycl(\pi, B^k)\times  B^m}=(Fix(\pi, D_k))^{B^m}
=D_m^{Cycl(\pi,B^k)}.$$
\begin{lem}
Suppose that a permutation $\pi$ is acting on $B^k$ and $n>k$.
Then when we consider $\pi$ acting on $B^n$, 
$$Fix(\pi,D_n)= (Fix(\pi, D_k))^{B^m} =D_m^{Cycl(\pi,B^k)}$$
where $m=n-k$.
\end{lem}
This lemma was used by Wiedemann~\cite{W} and by Fidytek, Mostowski, Somla and Szepietowski~\cite{FMSS} in order to compute $d_8=|D_8|$
and by Pawelski~\cite{P} in order to count and generate fixes of several permutations  acting on $D_n$.

\section{Applications}
In this section we present a few examples.
\subsection{$D_2$}
Consider the poset $D_2=B^{B^2}=\{0000,0001,0011,0101,0111,1111\}$ and its array
$$
M(D_2)=\left(
\begin{array}{cc|cc|cc}
1& 1& 1& 1& 1& 1 \\
0& 1& 1& 1& 1& 1 \\
\hline
0& 0& 1& 0& 1& 1 \\
0& 0& 0& 1& 1& 1 \\
\hline
0& 0& 0& 0& 1& 1 \\
0& 0& 0& 0& 0& 1
\end{array}
\right)
$$

Consider the array 
$$
M(D_2)^2=\left(
\begin{array}{cc|cc|cc}
1& 2& 3& 3& 5& 6 \\
0& 1& 2& 2& 4& 5 \\
\hline
0& 0& 1& 0& 2& 3 \\
0& 0& 0& 1& 2& 3 \\
\hline
0& 0& 0& 0& 1& 2 \\
0& 0& 0& 0& 0& 1
\end{array}
\right)
$$

$Sum(M(D_2)^2)=50$ which is equal to $|D_2^{P_3}|=|(B^{B^2})^{P_3}|=|B^{B^2\times P_3}|=|Fix((12),D_4)|$.

$SumSq(M(D_2)^2)=168$ which is equal to $|D_2^{B^2}|=|(B^{B^2})^{B^2}|  =|B^{B^2\times B^2}|=|B^{B^4}|=|D_4|=d4$.

$$
M(D_2)^3=\left(
\begin{array}{cc|cc|cc}
1& 3& 6& 6& 14& 20 \\
0& 1& 3& 3& 9& 14 \\
\hline
0& 0& 1& 0& 3& 6 \\
0& 0& 0& 1& 3& 6 \\
\hline
0& 0& 0& 0& 1& 3 \\
0& 0& 0& 0& 0& 1
\end{array}
\right)
$$

$Sum(M(D_2)^3)=105$ which is equal to 
$|D_2^{P_4}|=|(B^{B^2})^{P_4}|=|B^{B^2\times P_4}|=|F((123),D_5)|$.

\subsection{Permutation $(12)$}
Consider the permutation $\pi=(12)$. When $\pi$  acts on $B^2$, we have  
three cycles $Cycl((12),B^2)=P_3$ and $Fix((12),D_2)=B^{P_3}=P_4$.
When $(12)$ acts on $B^3$, then
$Cycl((12),B^3)=P_3\times B$ and $Fix((12),D_3)=B^{P_3\times B}$.
The nunber of fixes can be computed either by:
$$|Fix((12),D_3)|=|B^{P_3\times B}|=|(B^{P_3})^B|=|P_4^B|=Sum(M(P_4))=10,$$
or by
$$|Fix((12),D_3)|=|B^{P_3\times B}|=|(B^B)^{P_3}|=|D_1^{P_3}|=Sum(M(D_1)^2)=10.$$
In a similar way we can count fixes when $(12)$ acts on $B^4$, then
$Cycl((12),B^4)=P_3\times B^2$ and $Fix((12),D_4)=B^{P_3\times B^2}$.
The nunber of fixes can be computed either by:
$$|Fix((12),D_4)|=|B^{P_3\times B^2}|=|(B^{P_3})^{B^2}|=|P_4^{B^2}|=SumSq(M(P_4)^2)=50,$$
or by
$$|Fix((12),D_4)|=|B^{P_3\times B^2}|=|(B^{B^2})^{P_3}|=|D_2^{P_3}|=Sum(M(D_2)^2)=50.$$

\subsection{Permutation $(123)$}
Consider the permutation $\pi=(123)$. When $\pi$  acts on $B^3$, we have  
four cycles $Cycl((123),B^3)=P_4$ and $Fix((123),D_3)=B^{P_4}=P_5$.
When $(123)$ acts on $B^4$, then
$Cycl((123),B^4)=P_4\times B$ and $Fix((123),D_4)=B^{P_4\times B}$.
The nunber of fixes can be computed either by:
$$|Fix((123),D_4)|=|B^{P_4\times B}|=|(B^{P_4})^B|=|P_5^B|=Sum(M(P_5))=15,$$
or by
$$|Fix((123),D_4)|=|B^{P_4\times B}|=|(B^B)^{P_4}|=|D_1^{P_4}|=Sum(M(D_1)^3)=15.$$
In a similar way we can count fixes when $(123)$ acts on $B^5$, then
$Cycl((123),B^5)=P_4\times B^2$ and $Fix((123),D_5)=B^{P_4\times B^2}$.
The nunber of fixes can be computed either by:
$$|Fix((123),D_5)|=|B^{P_4\times B^2}|=|(B^{P_4})^{B^2}|=|P_5^{B^2}|=SumSq(M(P_5)^2)=105,$$
or by
$$|Fix((123),D_5)|=|B^{P_4\times B^2}|=|(B^{B^2})^{P_4}|=|D_2^{P_4}|=Sum(M(D_2)^3)=105.$$

\subsection{Permutation $(1234)$}
The permutation $(1234)$ acting on $B^4$ has six cycles:
$C_0=\{0000\}$,
$C_1=\{1000,0001,0010,0100\}$,
$C_2=\{1100,1001,0011,0110\}$,
$C_3=\{1010,0101\}$,
$C_4=\{1110,1101,1011,0111\}$,
and $C_5=\{1111\}$.
They are of length 1, 2, or 4 and ordered by
$C_0<C_1<C_2,C_3<C_4<C_5$.
There are 8 fixes of $(1234)$ acting on $D_4$ with the array
$$
M(Fix((1234),D_4))=\left(
\begin{array}{ccc|cc|ccc}
1& 1& 1& 1& 1& 1& 1& 1\\
0& 1& 1& 1& 1& 1& 1& 1 \\
0& 0& 1& 1& 1& 1& 1& 1 \\
\hline
0& 0& 0& 1& 0& 1& 1& 1 \\
0& 0& 0& 0& 1& 1& 1& 1 \\
\hline
0& 0& 0& 0& 0& 1& 1& 1 \\
0& 0& 0& 0& 0& 0& 1& 1 \\
0& 0& 0& 0& 0& 0& 0& 1 \\
\end{array}
\right)
$$
$Sum(M(Fix((1234),D_4)))=35=|Fix((1234),D_5))|$.

$$
M(Fix((1234),D_4))^2=\left(
\begin{array}{ccc|cc|ccc}
1& 2& 3& 4& 4& 6& 7& 8\\
0& 1& 2& 3& 3& 5& 6& 7 \\
0& 0& 1& 2& 2& 4& 5& 6 \\
\hline
0& 0& 0& 1& 0& 2& 3& 4 \\
0& 0& 0& 0& 1& 2& 3& 4 \\
\hline
0& 0& 0& 0& 0& 1& 2& 3 \\
0& 0& 0& 0& 0& 0& 1& 2 \\
0& 0& 0& 0& 0& 0& 0& 1 \\
\end{array}
\right)
$$
$SumSq(M(Fix((1234),D_4))^2)=494=|Fix((1234),D_6))|$.

$$
M(Fix((1234),D_4))^3=\left(
\begin{array}{ccc|cc|ccc}
1&3& 6& 10& 10& 20& 27& 35\\
0& 1& 3& 6& 6& 14& 20&27 \\
0& 0& 1& 3& 3& 9& 14& 20 \\
\hline
0& 0& 0& 1& 0& 3& 6& 10 \\
0& 0& 0& 0& 1& 3& 6& 10 \\
\hline
0& 0& 0& 0& 0& 1& 3& 6 \\
0& 0& 0& 0& 0& 0& 1& 3 \\
0& 0& 0& 0& 0& 0& 0& 1 \\
\end{array}
\right)
$$
$Sum(M(Fix((1234),D_4))^3)=294=|Fix((1234)(567),D_7))|$.

\subsection{Permutation $(12345)$}
The permutation $(12345)$ acting on $B^5$ has eight cycles:
$C_0=\{00000\}$,
$C_1=\{10000,00001,00010,00100,01000\}$,
$C_2=\{11000,10001,00011,00110,01100\}$,
$C_3=\{10100,01001,10010,00101,01010\}$,
$C_4=\{11100,11001,10011,00111,01110\}$,
$C_5=\{10110,01101,11010,10101,01011\}$,
$C_6=\{11110,11101,11011,10111,01111\}$,
and $C_7=\{11111\}$.
They are of length 1 or 5 and ordered by
$C_0<C_1<C_2,C_3<C_4,C_5<C_6<C_7$.
There are 11 fixes of $(12345)$ acting on $D_5$ with the array
$$
M(Fix((12345),D_5))=\left(
\begin{array}{ccc|cc|c|cc|ccc}
1& 1& 1& 1& 1& 1& 1& 1& 1& 1& 1\\
0& 1& 1& 1& 1& 1& 1& 1& 1& 1& 1\\
0& 0& 1& 1& 1& 1& 1& 1& 1& 1& 1\\
\hline
0& 0& 0& 1& 0& 1& 1& 1& 1& 1& 1\\
0& 0& 0& 0& 1& 1& 1& 1& 1& 1& 1\\
\hline
0& 0& 0& 0& 0& 1& 1& 1& 1& 1& 1\\
\hline
0& 0& 0& 0& 0& 0& 1& 0& 1& 1& 1\\
0& 0& 0& 0& 0& 0& 0& 1& 1& 1& 1\\
\hline
0& 0& 0& 0& 0& 0& 0& 0& 1& 1& 1\\
0& 0& 0& 0& 0& 0& 0& 0& 0& 1& 1\\
0& 0& 0& 0& 0& 0& 0& 0& 0& 0& 1\\
\end{array}
\right)
$$
$Sum(M(Fix((12345),D_5)))=64=|Fix((12345),D_6))|$.

$$
M(Fix((12345),D_5))^2=\left(
\begin{array}{ccc|cc|c|cc|ccc}
1& 2& 3& 4& 4& 6& 7& 7& 9&10& 11\\
0& 1& 2& 3& 3& 5& 6& 6& 8& 9& 10\\
0& 0& 1& 2& 2& 4& 5& 5& 7& 8& 9\\
\hline
0& 0& 0& 1& 0& 2& 3& 3& 5& 6& 7\\
0& 0& 0& 0& 1& 2& 3& 3& 5& 6& 7\\
\hline
0& 0& 0& 0& 0& 1& 2& 2& 4& 5& 6\\
\hline
0& 0& 0& 0& 0& 0& 1& 0& 2& 3& 4\\
0& 0& 0& 0& 0& 0& 0& 1& 2& 3& 4\\
\hline
0& 0& 0& 0& 0& 0& 0& 0& 1& 2& 3\\
0& 0& 0& 0& 0& 0& 0& 0& 0& 1& 2\\
0& 0& 0& 0& 0& 0& 0& 0& 0& 0& 1\\
\end{array}
\right)
$$
$Sum(M(Fix((12345),D_4))^2)=264=|Fix((12345)(67),D_7))|$

$SumSq(M(Fix((12345),D_4))^2)=1548=|Fix((12345),D_7))|$

$$
M(Fix((12345),D_5))^3=\left(
\begin{array}{ccc|cc|c|cc|ccc}
1& 3& 6&10&10&20&27&27&43&53&64\\
0& 1& 3& 6& 6&14&20&20&34&43&53\\
0& 0& 1& 3& 3& 9&14&14&26&34&43\\
\hline
0& 0& 0& 1& 0& 3& 6& 6&14&20&27\\
0& 0& 0& 0& 1& 3& 6& 6&14&20&27\\
\hline
0& 0& 0& 0& 0& 1& 3& 3& 9&14&20\\
\hline
0& 0& 0& 0& 0& 0& 1& 0& 3& 6&10\\
0& 0& 0& 0& 0& 0& 0& 1& 3& 6&10\\
\hline
0& 0& 0& 0& 0& 0& 0& 0& 1& 3& 6\\
0& 0& 0& 0& 0& 0& 0& 0& 0& 1& 3\\
0& 0& 0& 0& 0& 0& 0& 0& 0& 0& 1\\
\end{array}
\right)
$$
$Sum(M(Fix((12345),D_5))^3)=870=|Fix((12345)(678),D_8))|$.

\section{Permutation $(12)(34)$}
 Consider partition $A_4=\{1,2,3,4\}=\{1,2\}+\{3,4\}$,
and two permutations $\pi=(12)$ and $\rho=(34)$. Furthermore, consider partition of $B^4$ into four subcubes:

$B^4_{00}=\{0000,1000,0100,1100\}=B^2\times\{00\}$

$B^4_{10}=\{0010,1010,0110,1110\}=B^2\times\{10\}$

$B^4_{01}=\{0001,1001,0101,1101\}=B^2\times\{01\}$

$B^4_{11}=\{0011,1011,0111,1111\}=B^2\times\{11\}$

Each of these subcubes is isomorphic to $B^2$.
There are three kinds of cycles of $\pi\circ\rho$ acting on $B^4$:
\begin{itemize}
\item $(0000)$, $(1000,0100)$, $(1100)$. They are contained in $B^4_{00}$ and are 
isomorphic to the cycles of $\pi$ acting on $B^2$.
\item $(0011)$, $(1011,0111)$, $(1111)$. They are contained in $B^4_{11}$ and are isomorphic to the cycles of $\pi$ acting on $B^2$.
\item $(0010,0001)$, $(1010,0101)$, $(0110,1001)$, $(1110.1101)$.
Each of the cycles contains two elements $\{x,y\}$ such that $x\in B^4_{10}$, $y\in B^4_{01}$, and $y=\pi\circ\rho(x)$. Moreover each $x\in B^4_{10}$ belongs to one of these cycles.
\end{itemize}

Suppose that $f$ is a fix of $\pi\circ\rho$ acting on $D_4$ and consider four restrictions:

$f_{00}=f|_{B^4_{00}}$, 
$f_{10}=f|_{B^4_{10}}$,
$f_{01}=f|_{B^4_{01}}$,
$f_{11}=f|_{B^4_{11}}$.

They satisfy the following conditions:
\begin{enumerate}
\item $f_{00},f_{11}\in Fix(\pi,D_2).$ Here we identify
$B^2\times\{00\}$ (and $B^2\times\{11\}$) with $B^2$ and functions 
$B^{B^2\times\{00\}}$ (and $B^{B^2\times\{11\}}$) with $B^{B^2}$.
\item $f_{10},f_{01}\in B^{B^2}=D_2$
Here we identify  functions 
$B^{B^2\times\{10\}}$ (and $B^{B^2\times\{01\}}$) with $B^{B^2}$.
\item $f_{10}=\pi(f_{01})$
\item $f_{00}\le f_{10},f_{10}\le f_{11}$.
\end{enumerate}
On the other hand if for a function $f$ its restrictions 
$f_{00}=f|_{B^4_{00}}$, 
$f_{10}=f|_{B^4_{10}}$,
$f_{01}=f|_{B^4_{01}}$,
$f_{11}=f|_{B^4_{11}}$ satisfy conidions (1--4) then 
$f$ is a fix of $\pi\circ\rho$ acting on $B^4$.
\section{Permutation $(12)(34)(56)(78)$}
 Consider partition $A_4=\{1,\dots n+2\}=\{1,\dots n\}+\{n+1,n+2\}$,
and two permutations $\pi$ acting on $\{1,\dots n\}$ and $\rho=(n+1,n+2)$.
and suppose that cycles of $\pi$ are of length 1 or 2.
Consider partition of $B^{n+2}$ into four subcubes:

$B^{n+2}_{00}=B^n\times\{00\}$

$B^{n+2}_{10}=B^n\times\{10\}$

$B^{n+2}_{01}=B^n\times\{01\}$

$B^{n+2}_{11}=B^n\times\{11\}$

There are three kinds of cycles of $\pi\circ\rho$ acting on $B^{n+2}$:
\begin{itemize}
\item Those contained in $B^{n+2}_{00}$; isomorphic to the cycles of $\pi$ acting on $B^n$.
\item Those contained in $B^{n+2}_{11}$; isomorphic to the cycles of $\pi$ acting on $B^n$.
\item 
Each $x\in B^{n+2}_{10}$ belongs to the cycle with $y=\pi\circ\rho(x)\in B^{n+2}_{01}$.
\end{itemize}

Suppose that $f$ is a fix of $\pi\circ\rho$ acting on $D_{n+2}$ and consider 
four restrictions:

$f_{00}=f|_{B^{n+2}_{00}}$, 
$f_{10}=f|_{B^{n+2}_{10}}$,
$f_{01}=f|_{B^{n+2}_{01}}$,
$f_{11}=f|_{B^{n+2}_{11}}$.

They satisfy the following conditions:
\begin{enumerate}
\item $f_{00},f_{11}\in Fix(\pi,D_n)$
\item $f_{10},f_{01}\in B^{B^n}=D_n$
\item $f_{10}=\pi(f_{01})$
\item $f_{00}\le f_{10},f_{10}\le f_{11}$.
\end{enumerate}
On the other hand if for a function $f$ its restrictions
$f_{00}=f|_{B^{n+2}_{00}}$, 
$f_{10}=f|_{B^{n+2}_{10}}$,
$f_{01}=f|_{B^{n+2}_{01}}$,
$f_{11}=f|_{B^{n+2}_{11}}$ satisfy conidions (1--4) then $f$ is a fix of $\pi\circ\rho$ acting on $D_{n+2}$.

\noindent{\bf Algorithm counting fixes}\\
Input: posets $D_n=B^{B^n}$ and $Fix(\pi,D_n)$.\\
Output: $|Fix(\pi\circ\rho, D_{n+2})|$.

\begin{itemize}
\item $Sum := 0$
\item For each $f_{10}\in D_n$:
\begin{itemize}
\item $f_{01}:=\pi(f_{10})$;
\item $Down:= |\{g\in Fix(\pi,D_n): g\le f_{10}\&f_{01}\}|$\qquad//the number of possibilities for choosing $f_{00}$
\item $Up:= |\{g \in Fix(\pi,D_n): g\ge f_{10}|f_{01}\}|$ \qquad//the number of possibilities for choosing $f_{11}$
\item $Sum:= Sum + Down\cdot Up$
\end{itemize}
\item Return $|Fix(\pi\circ\rho, D_{n+2})|:=Sum$.
\end{itemize}
Note that for each function $g\in g\in Fix(\pi,D_n)$:
   
$g\le f_{10}$ and $g\le f_{01}$ if and only if $g\le f_{10}\&f_{01}$

and

$g\ge f_{10}$ and $g\ge f_{01}$ if and only if $g\ge f_{10}|f_{01}$.

Similar algorithm was used by Pawelski~\cite{P} in order to count  fixes of the permutation $(12)(34)(56)(78)$ acting
on $D_8$.

\begin{ex}
Consider agorithm working for the permutation $(12)(34)$ acting on $D_4$.

$D_2=\{0000<0001<0011,0101<0111<1111,\}$

$Fix((12),D_2)=\{0000<0001<0111<1111\}$
\begin{itemize}
\item for $f_{10}=0000$:  $f_{01}=0000$; 
$f_{10}\&f_{01}=0000$; $Down=1$; 
$f_{10}|f_{01}=0000$; $Up=4$.
\item for $f_{10}=0001$: $f_{01}=0001$;
$f_{10}\&f_{01}=0001$; $Down=2$;
$f_{10}|f_{01}=0001$; $Up=3$.
\item for $f_{10}=0011$: $f_{01}=0101$;
 $f_{10}\&f_{01}=0001$; $Down=2$;
 $f_{10}|f_{01}=0111$; $Up=2$.
\item for $f_{10}=0101$: $f_{01}=0011$; 
$f_{10}\&f_{01}=0001$; $Down=2$; 
$f_{10}|f_{01}=0111$; $Up=2$.
\item for $f_{10}=0111$: $f_{01}=0111$; 
$f_{10}\&f_{01}=0111$; $Down=3$;
 $f_{10}|f_{01}=0111$; $Up=2$.
\item for $f_{10}=1111$: $f_{01}=1111$; 
$f_{10}\&f_{01}=1111$; $Down=4$; 
$f_{10}|f_{01}=1111 $; $Up=1$.
\end{itemize}
Algorithm returns $|Fix((12)(34),D_4)|=
1\cdot 4+2\cdot 3+2\cdot 2+2\cdot 2+3\cdot 2+4\cdot 1=28$.
\end{ex}
 
\section{Generating $Fix$}
In this section we present one more method to generate  $Fix(\pi,D_n)$ fixes of a permutation $\pi$ acting on $D_n$.
We start with the poset $Cycl(\pi,B^n)$ with its array $M(Cycl(\pi,B^n))$
For example, consider the permutation $(12)$  acting on $B^3$. The poset
$Cycl((12),B^3)=\{a<b<c\}\times\{0<1\}=\{a0,b0,c0,a1,b1,c1\}$
has the matrix
$$
M(Cycl((12)),B^3))=\left(
\begin{array}{cccccc}
1& 1& 1& 1& 1& 1\\
0& 1& 1& 0& 1& 1\\
0& 0& 1& 0& 0& 1\\
1& 0& 0& 1& 1& 1\\
0& 0& 0& 0& 1& 1\\
0& 0& 0& 0& 0& 1\\
\end{array}
\right)
$$
We shall identify rows of the array with subsets of $Cycl(\pi,B^n)$
and with functions from $Cycl(\pi,B^n)$ to $ \{0,1\}$.
It is well known that monotone functions from $Cycl(\pi,B^n)$ to $\{0,1\}$ may be identified with upsets.
A subset $U\subset Cycl(\pi,B^n)$ is an upset if 

for every $x,y$,  $x\in U$ and $x\le y$ implies $y\in U$.\\

\noindent
Each row in the array $M(Cycl(\pi,B^n))$ represents the upset
$Up(c)=\{x\in C:x\ge c\}$.
The set of all upsets can be generated in the following way:
We start with rows of the array $M(Cycl(\pi,B^n))$. Then we add the zero vector
and  or $x|y$ of every pair $x,y$ already in $Fix$.

\bigbreak
\noindent Algorithm generating $Fix(\pi,B^n)$\\
Input: poset $C=Cycl(\pi, B^n)$ and its array\\
Output: $Fix(\pi,D_n)$.
\begin{itemize}
\item $Fix:=\emptyset$
\item add zero vector to $Fix$
\item For each $c\in C$:
\begin{itemize}
\item for each $x\in Fix$ add $x|Up(c)$ to $Fix$
\item remove repetitions in $Fix$
\end{itemize}
\item Return $Fix(\pi,D_{n}):=Fix$
\end{itemize}
For example, the algorithm adds four rows to
the array $M(Cycl((12),B^3))$
$$
\begin{array}{cccccc}
1& 1& 1& 1& 1& 1\\
0& 1& 1& 0& 1& 1\\
0& 0& 1& 0& 0& 1\\
1& 0& 0& 1& 1& 1\\
0& 0& 0& 0& 1& 1\\
0& 0& 0& 0& 0& 1\\
\hline
\hline
0& 0& 0& 0& 0& 0\\
0& 1& 1& 1& 1& 1\\
0& 0& 1& 1& 1& 1\\
0& 0& 1& 0& 1& 1\\
\end{array}
$$
These ten rows form the poset $Fix((12),D_3)$ with the  partial order
defined by
$$x\le y\hbox{ iff } x|y=y.$$

\section{Fixes}
In this section we present tables with numbers of fixes of all permutations acting in $D_n$ for $n=3,...,8$.
Values for $n\le 6$ are from~\cite{CS1}, values for $n=7,8$  are from~\cite{CS2,P}.

$n=3$
\begin{tabular}[h]{|c|c|c|c|}
\hline
 $i$       &$\pi_i$&$\mu_i$ &$|Fix(\pi_i,D_3)|$\\
\hline
1& $e$      & 1 & 20\\
2& (12)     & 3 & 10 \\
3& (123)    & 2 & 5 \\
\hline
\end{tabular}

\bigbreak
$n=4$
\begin{tabular}[h]{|c|c|c|c|}
\hline
 $i$       &$\pi_i$&$\mu_i$ &$|Fix(\pi_i,D_4)|$\\
\hline
1& $e$      & 1  & 168\\
2& (12)     & 6  & 50  \\
3& (123)    & 8 & 15 \\ 
4&  (1234)  &  6 &  8 \\
5&  (12)(34)&  3 & 28\\
\hline
\end{tabular}

\bigbreak
$n=5$
\begin{tabular}[h]{|c|c|c|c|}
\hline
 $i$       &$\pi_i$&$\mu_i$ &$|Fix(\pi_i,D_5)|$\\
\hline
1& $e$      & 1  & 7 581\\
2& (12)     & 10  & 887  \\
3& (123)    & 20 & 105 \\ 
4&  (1234)  &  30 &  35 \\
5&  (12)(34)&  15 & 309\\
6& (12345)  & 24& 11\\
7& (12)(345)& 20 &35\\
\hline
\end{tabular}

\bigbreak
$n=6$
\begin{tabular}[h]{|c|c|c|c|}
\hline
 $i$       &$\pi_i$&$\mu_i$ &$|Fix(\pi_i,D_6)|$\\
\hline
1& $e$         &1   & 7 828 354\\
2& (12)        &15  & 160 948  \\
3& (123)       &40  & 3 490 \\ 
4& (1234)     &90  &  494 \\
5& (12)(34)   &45  & 24 302\\
6& (12345)     &144 & 64\\
7& (123456)    &120 & 44\\
8& (12)(345)   &120 & 490\\
9& (123)(456)  &40  & 562\\
10&(12)(3456)  &90  & 324\\
11&(12)(34)(56)&15  & 860\\
\hline
\end{tabular}

\bigbreak
$n=7$
\begin{tabular}[h]{|c|c|c|c|}
\hline
 $i$       &$\pi_i$&$\mu_i$ &$|Fix(\pi_i,D_7)|$\\
\hline
1& $e$            &1   & 2 414 682 040 998\\
2& (12)           &15  & 2 208 001 624  \\
3& (123)          &40  & 2 068 224 \\ 
4&  (1234)        &90  & 60 312 \\
5& (12345)        &144 & 1 548\\
6& (123456)       &120 & 766\\
7& (1234567)      &120 & 101\\
8&  (12)(34)      &45  & 67 922 470\\
9&  (12)(345)     &45  & 59 542\\
10& (12)(3456)    &120 & 26 878\\
11& (12)(34567)   &120 & 264\\
12& (123)(456)    &120 & 69 264\\
13& (123)(4567)   &120 & 294\\
14&(12)(34)(56)   &15  & 12 015 832 860\\
15&(12)(34)(567)  &15  & 10 192\\
\hline
\end{tabular}

\bigbreak
$n=8$
\begin{tabular}[h]{|c|c|c|c|}
\hline
 $i$       &$\pi_i$&$\mu_i$ &$|Fix(\pi_i,D_8)|$\\
\hline
1& $e$              &1    & 56 130 437 228 687 557 907 788\\
2& (12)             &28   & 101 627 867 809 333 596  \\
3& (123)            &112  & 262 808 891 710 \\ 
4& (1234)           &420  & 424 234 996 \\
5& (12345)          &1344 & 531 708\\
6& (123456)         &3366 & 144 320\\
7& (1234567)        &5760 & 3 858\\
8& (12345678)       &5040 & 2 364\\
9&  (12)(34)        &210 & 182 755 441 509 724\\
10& (12)(345)       &1120& 401 622 018\\
11& (12)(3456)      &2520& 93 994 196\\
12& (12)(34567)     &4032& 21 216\\
13& (12)(345678)    &3360& 70 096\\
14& (123)(456)      &1120& 535 426 780\\
15& (123)(4567)     &3360& 25 168\\
16& (123)(45678)    &2688& 870\\
17& (1234)(5678)    &1260& 3 211 276\\
18&(12)(34)(56)     &420 & 7 377 670 895 900\\
19&(12)(34)(567)    &1680& 16 380 370\\
20&(12)(34)(5678)   &1260& 37 834 164\\
21&(12)(345)(678)   &1120& 3 607 596\\
22&(12)(34)(56)(78) &105 & 2 038 188 253 420\\
\hline
\end{tabular}

\section{Known values of $d_n$ and $r_n$}

\begin{tabular}[h]{|c|c|c|}
\hline
 $n$       &$d_n$&$r_n$\\
\hline
0&  2       &   2\\  
1&  3       &   3\\
2&  6       &   5\\  
3&  20       &   10\\  
4&  168       &   30\\
5&  7 581       &   210\\  
6&  7 828 354   &   16 353\\  
7&  2 414 682 040 998  & 490 013 148  \\
8&  56 130 437 228 687 557 907 788 &  1 392 195 548 889 993 358\\
\hline
\end{tabular}


\end{document}